\documentclass[12pt]{article}
\usepackage{amsmath, amsfonts}
\numberwithin{equation}{section}
\begin{document}

\author{Ajai Choudhry}
\title{Some diophantine problems concerning\\ a pair of rational triangles\\ with a common circumradius}
\date{}
\maketitle

\begin{abstract}
A triangle with rational sides and rational area is
called a rational triangle. In this paper we consider three problems of finding pairs of rational triangles which have a common circumradius as well as either a common perimeter or a common inradius or a common area. While several similar problems concerning pairs of rational triangles have been considered in the existing literature, these three problems have not been studied till now. For each of these problems, we give a parametric solution and we also indicate how additional parametric solutions of these problems may be obtained.
\end{abstract}

Mathematics Subject Classification 2020: 11D25

Keywords: rational triangles; Heron triangles; circumradius of a triangle; inradius of a triangle. 

\section{Introduction}

\hspace{0.25in} A triangle with rational sides and rational area is
called a rational triangle. Diophantine problems concerning rational triangles have attracted considerable attention.
For instance, several mathematicians have considered the problem of finding two rational triangles with  a common   perimeter and a common area (see \cite{Aa}, \cite{Br}, \cite{HM}, \cite{Lu}, \cite{Yi}). In fact, Choudhry \cite{Ch} has described a method of generating an arbitrarily large number of scalene rational triangles with a common   perimeter and a common area.  Skalba and Ulas have considered the problem of  finding pairs of Pythagorean triangles with given ratios between catheti. Regarding problems concerning rational triangles with a common circumradius, it has been shown by Lehmer \cite[Theorem XI, p. 101]{Le} that there exist infinitely many rational triangles with a common circumradius.  Further,  Andrica and \c{T}urca\c{s} \cite{AT} have recently proved that there  are no pairs consisting of a rational right triangle and a rational isosceles triangle which have the same  circumradius and the same inradius or which have the  same circumradius and the same perimeter. 

This paper is concerned with three diophantine problems pertaining to a pair of scalene rational triangles that have the same circumradius. The three  problems require that we find pairs of rational triangles with:

\noindent (i) a common circumradius and a common perimeter;

 \noindent (ii) a common circumradius and a common inradius;

\noindent (iii) a common circumradius and a common area.

We note that none of the above three  problems has been considered earlier in the literature. We obtain parametric solutions of each of the above problems. We also show how more parametric solutions of these problems may be obtained. We note that rational values of the parameters may  yield two triangles with rational sides as a solution to any of our  three diophantine  problems. In each case, we may, after appropriate scaling, readily obtain two triangles  whose  sides and  areas are given by integers and which have the desired properties.

\section{Some basic formulae regarding rational triangles}\label{formulae}
In this section we will give basic formulae  for the sides, the area, the circumradius and the inradius of a general triangle whose sides and area are rational. 

We note that Brahmagupta \cite[p. 191]{Di} and Euler \cite[p. 193]{Di} have independently given two sets of formulae for the sides of a general rational triangle. The problem of determining all rational triangles has also been considered by Carmichael \cite[pp. 11--13]{Ca} and by Lehmer \cite{Le}. We could try to use these well-known formulae about rational triangles to find pairs of such triangles with  a common  circumradius as well as a common  perimeter/ common inradius. The resulting diophantine equations are, however,  difficult to solve.

We will now derive a new set of formulae for the sides and area of an arbitrary  rational triangle.  It is interesting to observe that using  the new formulae given below we can neatly resolve the problems of finding pairs of rational triangles with a common  circumradius as well as a common  perimeter/ common inradius.

Let  $a, b , c$, be  the sides of an arbitrary  rational triangle.  The area, circumradius and inradius of the  triangle, denoted by $A, R$ and $r$, respectively, are given by the following well-known formulae:
\begin{align}
A&=\sqrt{(a+b+c)(a+b-c)(b+c-a)(c+a-b)}/4, \label{defA}\\
R&=abc/\sqrt{(a+b+c)(a+b-c)(b+c-a)(c+a-b)}, \label{defR}\\
r&=\sqrt{(a+b+c)(a+b-c)(b+c-a)(c+a-b)}/\{2(a+b+c)\}. \label{defr}
\end{align}

On making  the invertible linear transformation defined by,
\begin{equation}
a=y+z,\quad b=z+x, \quad c=x+y, \label{ltt1}
\end{equation}
 the above formulae may be written as,

\begin{align}
A&=\sqrt{(x + y + z)xyz}, \label{defA2}\\
R&=(x + y)(y + z)(x + z)/\{4\sqrt{(x + y + z)xyz}\}, \label{defR2}\\
r&=\sqrt{(x + y + z)xyz}/(x + y + z), \label{defr2}
\end{align}
where we note that $x, y$ and $z$ are necessarily nonzero rational numbers.

It follows from \eqref{defA2} that the area of the triangle will be rational if and only if there is   a nonzero rational number $t$ such that 
\begin{equation}
(x + y + z)x=t^2yz, \label{relt}
\end{equation}
so that
\begin{equation}
z= (x + y)x/(t^2y - x), \label{valz}
\end{equation}
and now the values of $A, R$ and $r$ may be written, in terms of nonzero rational numbers $x, y$ and $t$, as follows:
\begin{align}
A&=txy(x + y)/(t^2y - x), \label{defA3}\\
R&=(x^2+t^2y^2)(t^2 + 1)/(4(t^2y - x)t), \label{defR3}\\
r&=x/t. \label{defr3}
\end{align}

We note that, using the relations \eqref{ltt1} and the value of $z$ given by \eqref{valz},  the sides $a, b, c$ of our triangle may be written, in terms of three arbitrary nonzero parameters $x, y$ and $t$, as follows:
\begin{equation}
a = (x^2+t^2y^2)/(t^2y - x),\quad  b = xy(t^2 + 1)/(t^2y - x),\quad  c = x + y. \label{valabcgen}
\end{equation}
The perimeter $P$ of the triangle is now given, in terms of $x, y$ and $t$, by the formula,
\begin{equation}
P=2t^2y(x + y)/(t^2y - x). \label{defP}
\end{equation}

To find pairs of triangles with the desired properties, we will begin with two triangles whose sides may be written, using the formulae \eqref{valabcgen},  in terms of arbitrary parameters $x_i, y_i, t_i,\,i=1,2$. We will then impose the desired conditions on the two triangles, and solve the resulting diophantine equations. We will follow this approach in Sections \ref{eqRs} and \ref{eqRr} to obtain pairs of rational triangles with a common circumradius and a common perimeter/  common inradius.

\section{Three diophantine problems concerning rational triangles with the same circumradius}\label{diophprob}

\subsection{Pairs of rational triangles with a common  circumradius and  a common perimeter}\label{eqRs}

We will now obtain  examples of pairs of triangles with a common  circumradius and  a common perimeter.  

Let the sides $a_1, b_1, c_1$ and $a_2, b_2, c_2$ of the two triangles, respectively, be expressed, using the formulae \eqref{valabcgen}, in terms of two sets of arbitrary nonzero rational parameters $x_i, y_i, t_i, \;i=1,2$, applicable to the two triangles, respectively.

It follows from the formula \eqref{defR3} that if the two triangles have a common circumradius, the parameters  $x_i, y_i, t_i, \; i=1, 2$, must satisfy the following condition:
\begin{equation}
(x_1^2+t_1^2y_1^2)(t_1^2 + 1)/(4(t_1^2y_1 - x_1)t_1)=(x_2^2+t_2^2y_2^2)(t_2^2 + 1)/(4(t_2^2y_2 - x_2)t_2). \label{condR}
\end{equation}

Further,  on using the formula \eqref{defP}, the condition that the two  triangles also have the same perimeter may be be written as follows:
\begin{equation}
2t_1^2y_1(x_1 + y_1)/(t_1^2y_1 - x_1) = 2t_2^2y_2(x_2 + y_2)/(t_2^2y_2 - x_2). \label{condP}
\end{equation}

We will now solve the simultaneous equations \eqref{condR} and \eqref{condP}. On equating either side of Eq.~\eqref{condP} to $m$, and solving for $x_1$ and $x_2$, we get,
\begin{equation}
x_1 = t_1^2y_1(m - 2y_1)/(2t_1^2y_1 + m),\quad x_2 = t_2^2y_2(m - 2y_2)/(2t_2^2y_2 + m). \label{valx12Rp}
\end{equation}

On substituting the values of $x_1$ and $x_2$ given by \eqref{valx12Rp}, Eq.~\eqref{condR} may be written as follows:
\begin{equation}
(t_1^2 + 1)(4y_1^2t_1^2 + m^2)(2y_2t_2^2 + m)t_2 = (t_2^2 + 1)(4y_2^2t_2^2 + m^2)(2y_1t_1^2 + m)t_1. \label{condRp1}
\end{equation}
Eq.~\eqref{condRp1} may be considered as a cubic curve in $y_1$ and $y_2$, and a rational point on it, obtained by equating each side of \eqref{condRp1} to 0, is $(y_1, y_2)=(-m/(2t_1^2), -m/(2t_2^2))$. By drawing a tangent to the curve at this point, and taking its intersection with the curve \eqref{condRp1}, we get a rational solution of  \eqref{condRp1}. 

We can now obtain a solution of the simultaneous equations \eqref{condR} and \eqref{condP}, and using the relations \eqref{valabcgen}, we get the sides of the two triangles which have a common  circumradius and common perimeter. We omit the tedious details, and simply give below the sides $a_1, b_1, c_1$, and $a_2, b_2, c_2$ of the two triangles, obtained after appropriate scaling:
\begin{equation}
\begin{aligned}
a_1& =t_1(t_2^2 + 1)(t_1^4t_2^4 + 3t_1^4t_2^2 + 4t_1^3t_2^3 + 3t_1^2t_2^4 + t_1^4+ 2t_1^3t_2 + 3t_1^2t_2^2\\
 & \quad \quad  + 2t_1t_2^3 + t_2^4)(4t_1^6t_2^6 + 5t_1^6t_2^4 - 2t_1^5t_2^5+ 5t_1^4t_2^6 + 3t_1^6t_2^2\\
 & \quad \quad  - 2t_1^5t_2^3 + 2t_1^4t_2^4 - 2t_1^3t_2^5 + 3t_1^2t_2^6 + t_1^6 - 2t_1^3t_2^3 + t_2^6),\\
 b_1& = 2t_1t_2^3(t_1^2 + 1)^2(t_2^2 + 1)(t_1^2t_2^2 + t_1^2 + t_1t_2 + t_2^2)\\
 & \quad \quad \times (3t_1^5t_2^4 + t_1^4t_2^5 + 3t_1^5t_2^2 + t_1^4t_2^3 + t_1^3t_2^4 - t_1^2t_2^5\\
  & \quad \quad+ t_1^5 + t_1^4t_2 + t_1^3t_2^2 - t_1^2t_2^3 - t_1t_2^4 - t_2^5), \\
c_1 & = t_1(t_1 + t_2)(t_2^2 + 1)(3t_1^4t_2^4 + 3t_1^4t_2^2 + 3t_1^2t_2^4 + t_1^4+ t_1^2t_2^2 + t_2^4)\\
  & \quad \quad \times(2t_1^6t_2^5 + 2t_1^6t_2^3 - t_1^5t_2^4 + 3t_1^4t_2^5 - 3t_1^5t_2^2+ t_1^4t_2^3 \\
  & \quad \quad + t_1^3t_2^4 + 3t_1^2t_2^5 - t_1^5 - t_1^4t_2 - t_1^3t_2^2 + t_1^2t_2^3 + t_1t_2^4 + t_2^5),
\end{aligned}
\label{sidestrg1eqRp}
\end{equation} 
\begin{equation}
\begin{aligned}
a_2& =t_2(t_1^2 + 1)(t_1^4t_2^4 + 3t_1^4t_2^2 + 4t_1^3t_2^3 + 3t_1^2t_2^4 + t_1^4 + 2t_1^3t_2\\
 & \quad \quad + 3t_1^2t_2^2 + 2t_1t_2^3 + t_2^4)(4t_1^6t_2^6 + 5t_1^6t_2^4 - 2t_1^5t_2^5 + 5t_1^4t_2^6 \\
& \quad \quad + 3t_1^6t_2^2 - 2t_1^5t_2^3 + 2t_1^4t_2^4 - 2t_1^3t_2^5 + 3t_1^2t_2^6 + t_1^6 - 2t_1^3t_2^3 + t_2^6),\\
 b_2 &= 2t_1^3t_2(t_1^2 + 1)(t_2^2 + 1)^2(t_1^2t_2^2 + t_1^2 + t_1t_2 + t_2^2)\\
& \quad \quad \times (t_1^5t_2^4 + 3t_1^4t_2^5 - t_1^5t_2^2 + t_1^4t_2^3 + t_1^3t_2^4 + 3t_1^2t_2^5\\
& \quad \quad - t_1^5 - t_1^4t_2 - t_1^3t_2^2 + t_1^2t_2^3 + t_1t_2^4 + t_2^5),\\
c_2& =t_2(t_1 + t_2)(t_1^2 + 1)(3t_1^4t_2^4 + 3t_1^4t_2^2 + 3t_1^2t_2^4 + t_1^4 + t_1^2t_2^2 + t_2^4)\\
& \quad \quad \times (2t_1^5t_2^6 + 3t_1^5t_2^4 - t_1^4t_2^5 + 2t_1^3t_2^6 + 3t_1^5t_2^2 + t_1^4t_2^3 \\
& \quad \quad + t_1^3t_2^4 - 3t_1^2t_2^5 + t_1^5 + t_1^4t_2 + t_1^3t_2^2 - t_1^2t_2^3 - t_1t_2^4 - t_2^5).
\end{aligned}
\label{sidestrg2eqRp}
\end{equation} 

The common circumradius of the above two triangles is
\begin{multline}
\{(t_1^2 + 1)(t_2^2 + 1)(t_1^4t_2^4 + 3t_1^4t_2^2 + 4t_1^3t_2^3 + 3t_1^2t_2^4 + t_1^4 + 2t_1^3t_2 \\
+ 3t_1^2t_2^2 + 2t_1t_2^3 + t_2^4)(4t_1^6t_2^6 + 5t_1^6t_2^4 - 2t_1^5t_2^5 + 5t_1^4t_2^6 + 3t_1^6t_2^2\\
 - 2t_1^5t_2^3 + 2t_1^4t_2^4 - 2t_1^3t_2^5 + 3t_1^2t_2^6 + t_1^6 - 2t_1^3t_2^3 + t_2^6)\}/4,
\end{multline}
and the common perimeter is
\begin{multline}
4t_1^3t_2^3(t_1 + t_2)(t_1^2 + 1)(t_2^2 + 1)(t_1^2t_2^2 + t_1^2 + t_1t_2 + t_2^2)\\
\times (3t_1^4t_2^4 + 3t_1^4t_2^2 + 3t_1^2t_2^4 + t_1^4 + t_1^2t_2^2 + t_2^4). \quad \quad \quad \quad 
\end{multline}

As a numerical example, when $t_1=2, t_2=3$, we get, after appropriate scaling, two triangles with sides $ 1321940, 1166616,  1636180$ and $991455$,  $1548096, 1585185 $, having a common circumradius $1652425/2$ and a common perimeter $4124736$.

It is interesting to observe that when we take $t_1=1$, the first triangle becomes a right triangle while  taking $t_2=1$ makes the second  triangle  a right triangle. We give below the sides of the two triangles with $t_2=1$:
\begin{equation}
\begin{aligned}
a_1& = 2t_1(5t_1^4 + 6t_1^3 + 6t_1^2 + 2t_1 + 1)(13t_1^6 - 4t_1^5 + 7t_1^4 - 4t_1^3 + 3t_1^2 + 1),\\
 b_1& = 4t_1(t_1^2 + 1)^2(2t_1^2 + t_1 + 1)(7t_1^5 + 3t_1^4 + 2t_1^3 - 2t_1^2 - t_1 - 1),\\
 c_1& = 2t_1(t_1 + 1)(7t_1^4 + 4t_1^2 + 1)(4t_1^6 - 5t_1^5 + 3t_1^4 + 4t_1^2 + t_1 + 1), 
\end{aligned}
\label{sidestrg1eqRpspl}
\end{equation} 
\begin{equation}
\begin{aligned}
a_2& = (t_1^2 + 1)(5t_1^4 + 6t_1^3 + 6t_1^2 + 2t_1 + 1)(13t_1^6 - 4t_1^5 + 7t_1^4 \\
& \quad \quad - 4t_1^3 + 3t_1^2 + 1), \\
b_2& = -8t_1^3(t_1^2 + 1)(2t_1^2 + t_1 + 1)(t_1^5 - 3t_1^4 - 4t_1^2 - t_1 - 1),\\
 c_2& = (t_1 + 1)(t_1^2 + 1)(7t_1^4 + 4t_1^2 + 1)(9t_1^5 + t_1^4 + 4t_1^3 - 4t_1^2 - t_1 - 1).
\end{aligned}
\label{sidestrg2eqRpspl}
\end{equation} 

We now have a scalene triangle whose sides are given by \eqref{sidestrg1eqRpspl} and a right triangle whose sides are given by \eqref{sidestrg2eqRpspl} such that the two triangles have a common circumradius
\[(t_1^2 + 1)(5t_1^4 + 6t_1^3 + 6t_1^2 + 2t_1 + 1)(13t_1^6 - 4t_1^5 + 7t_1^4 - 4t_1^3 + 3t_1^2 + 1)/2\]
and a common perimeter
\[ 8t_1^3(t_1 + 1)(t_1^2 + 1)(2t_1^2 + t_1 + 1)(7t_1^4 + 4t_1^2 + 1).\]

 As a numerical example, when $t_1=2$, we get two triangles with sides $ 500516,  609400,  252324$ and   $625645,  123200,  613395$  having common circumradius $625645/2$ and common perimeter $1362240$.

We note that the cubic curve in $y_1$ and $y_2$  defined by Eq.~\eqref{condRp1} may be considered as an elliptic curve, and more rational points on the curve \eqref{condRp1} may be found by the well-known tangent and chord process or equivalently, by using the group law. These rational points will yield additional parametric solutions of the problem of finding pairs of rational triangles with a common circumradius and a common perimeter. 

\subsection{Pairs of rational triangles with a common circumradius and a common inradius}\label{eqRr}
We will now obtain two triangles with a common  circumradius and a common inradius. As in Section \ref{eqRs}, let the sides $a_i, b_i,  c_i, \, i=1, 2$, of the two triangles be expressed, using the formulae \eqref{valabcgen}, in terms of arbitrary parameters $x_i, y_i, t_i, \, i=1, 2$. 

The condition for the two triangles to have  a common circumradius is given by \eqref{condR} while, on using the formula \eqref{defr3}, the condition that they   have a common inradius may be written as follows:
\begin{equation}
x_1/t_1=x_2/t_2. \label{condr}
\end{equation}

We may therefore write $x_1=mt_1, x_2=mt_2$, where $m$ is an arbitrary parameter, and now the condition \eqref{condR} reduces to
\begin{multline}
t_2(t_1^2 + 1)y_1^2y_2 - t_1(t_2^2 + 1)y_1y_2^2 - m(t_1^2 + 1)y_1^2 + m(t_2^2 + 1)y_2^2\\
 - m^2t_1(t_2^2 + 1)y_1 + m^2t_2(t_1^2 + 1)y_2 - m^3(t_1 - t_2)(t_1 + t_2)=0. \label{condRr1}
\end{multline}

This is a cubic equation in $y_1$ and $y_2$, and it is easily observed that a rational point on the cubic curve defined by \eqref{condRr1} is given by $(y_1, y_2)=(t_2m, t_1m)$. We now draw a tangent to the cubic curve at the point $(t_2m, t_1m)$, and take its intersection with the curve, and thus obtain a new rational point on the curve \eqref{condRr1}. Using this new rational point, we readily obtain a rational solution of the simultaneous diophantine equations \eqref{condR} and \eqref{condr}. Now on using the relations \eqref{valabcgen}, we obtain the sides of two triangles with a common circumradius and a common inradius. On appropriate scaling, the sides $a_1, b_1, c_1$, and $a_2, b_2, c_2$ of the two triangles may be written as follows:
\begin{equation}
\begin{aligned}
a_1& = t_1(t_2^2 + 1)(t_1^2t_2^2 + t_1^2 - 8t_1t_2 + t_2^2 + 9),\\
 b_1 &= -(t_1t_2^2 - t_1 - 2t_2)(t_1^2 + 1)(2t_1t_2 - t_2^2 - 3), \\
c_1& = -2(t_1^2t_2^2 - t_1^2 - 4t_1t_2 + t_2^2 + 3)(t_1^2t_2 - 2t_1 - t_2),\\
\end{aligned}
\label{sidestrg1eqRr}
\end{equation} 
\begin{equation}
\begin{aligned}
a_2& = t_2(t_1^2 + 1)(t_1^2t_2^2 + t_1^2 - 8t_1t_2 + t_2^2 + 9),\\
 b_2& = (t_1^2t_2 - 2t_1 - t_2)(t_2^2 + 1)(t_1^2 - 2t_1t_2 + 3),\\
 c_2& = -2(t_1^2t_2^2 + t_1^2 - 4t_1t_2 - t_2^2 + 3)(t_1t_2^2 - t_1 - 2t_2),
\end{aligned}
\label{sidestrg2eqRr}
\end{equation} 
where $t_1$ and $t_2$ are arbitrary parameters.

The common circumradius of the above two triangles is 
\[(t_1^2 + 1)(t_2^2 + 1)(t_1^2t_2^2 + t_1^2 - 8t_1t_2 + t_2^2 + 9)/4,\]
while their common inradius is $2(t_1t_2^2 - t_1 - 2t_2)(t_1^2t_2 - 2t_1 - t_2)$.

As a numerical example, when $ t_1=9/2$ and $t2=7/6$, 
 we get, after appropriate scaling,  two triangles with sides $2055,  1105,  3002$, and $ 4795, 4845,  482$, with their common circumradius and common inradius being $58225/24$ and $228$, respectively.

As in the case of the two triangles in Section \ref{eqRs} given by \eqref{sidestrg1eqRp} and \eqref{sidestrg2eqRp}, we note that if we take $t_1=1$, the first triangle given by \eqref{sidestrg1eqRr} becomes a right triangle while on taking $t_2=1$, the second triangle given by \eqref{sidestrg2eqRr} becomes a right triangle. We give below the sides of the two triangles with $t_2=1$:
\begin{equation}
\begin{aligned}
a_1& = 2t_1(t_1^2 - 4t_1 + 5),\\ 
b_1 &= 2(t_1^2 + 1)(t_1 - 2),\\
 c_1& = 4(t_1 - 1)(t_1^2 - 2t_1 - 1)
\end{aligned}
\label{sidestrg1eqRrspl}
\end{equation} 
\begin{equation}
\begin{aligned}
a_2& = (t_1^2 + 1)(t_1^2 - 4t_1 + 5),\\
 b_2& = (t_1^2 - 2t_1 - 1)(t_1^2 - 2t_1 + 3), \\
c_2& = 4(t_1 - 1)^2.
\end{aligned}
\label{sidestrg2eqRrspl}
\end{equation}

We now have a scalene triangle whose sides are given by \eqref{sidestrg1eqRrspl} and a right triangle whose sides are given by \eqref{sidestrg2eqRrspl} such that the two triangles have a common circumradius
$(t_1^2 + 1)(t_1^2 - 4t_1 + 5)/2$ and a common inradius $2(t_1^2 - 2t_1 - 1)$.

As a numerical example, when $t_1=4$, 
 we get two triangles with sides $40, 68,  84$, and $85,  77,  36$, with the common circumradius and common inradius being $85/2$ and $14$, respectively.

As in Section \ref{eqRs}, the cubic curve defined by Eq.~\eqref{condRr1} may be considered as an elliptic curve, and additional rational points on this curve, found using the group law, will lead to additional parametric solutions of our problem. 
 
\subsection{Pairs of rational triangles with a common circumradius and a common area}\label{eqRA}
We will now find pairs of rational triangles with a common circumradius and a common area. It follows from formulae \eqref{defA} and \eqref{defR} that two  triangles, with rational  sides  $a_1, b_1, c_1$, and $a_2, b_2, c_2$,  respectively, will have a common circumradius and a common area if the following conditions are satisfied:
\begin{equation}
a_1b_1c_1=a_2b_2c_2, \label{cond1}
\end{equation}
and
\begin{multline}
(a_1+b_1+c_1)(a_1+b_1-c_1)(b_1+c_1-a_1)(c_1+a_1-b_1)\\
=(a_2+b_2+c_2)(a_2+b_2-c_2)(b_2+c_2-a_2)(c_2+a_2-b_2). \label{cond2}
\end{multline}
Further, the common area of the two triangles will be rational if and only if each side of Eq.~\eqref{cond2} is a perfect square.

To solve the simultaneous diophantine equations \eqref{cond1} and \eqref{cond2}, we write,
\begin{equation}
a_1=pu,\quad b_1=qv,\quad a_2=pv, \quad b_2=qu, \quad c_2=c_1, \label{subsabc}
\end{equation}
where $p, q, u $ and $ v$ are arbitrary parameters.  Now Eq.~\eqref{cond1} is identically satisfied while Eq.~\eqref{cond2} reduces to
\begin{equation}
(u - v)(u + v)(p - q)(p + q)\{2c_1^2 - (u^2 + v^2)(p^2 + q^2)\}=0. \label{cond2a}
\end{equation}

To obtain a nontrivial solution of Eqs.~\eqref{cond1} and \eqref{cond2}, the last factor on the left-hand side of Eq.~\eqref{cond2a} must be equated to 0. We now have   a quadratic equation in $u, v$ and $c_1$, and accordingly, we readily obtain the following solution of Eq.~\eqref{cond2a}:
\begin{equation}
\begin{aligned}
 u& = (m^2 + 2mn - n^2)p - (m^2 - 2mn - n^2)q,\\
 v& = (m^2 - 2mn - n^2)p + (m^2 + 2mn - n^2)q,\\
c_1&=(p^2 + q^2)(m^2 + n^2),
\end{aligned}
\label{valuvc1}
\end{equation}
where $m$ and  $n$ are arbitrary parameters.

On substituting the values of $u, v$ and $c_1$ in the relations \eqref{subsabc}, we obtain the following solution of  the simultaneous  equations \eqref{cond1} and \eqref{cond2}:
\begin{equation}
 \begin{aligned}
a_1 &= (m^2 + 2mn - n^2)p^2 - (m^2 - 2mn - n^2)pq,\\
 a_2& = (m^2 - 2mn - n^2)p^2 + (m^2 + 2mn - n^2)pq, \\
b_1& = (m^2 - 2mn - n^2)pq + (m^2 + 2mn - n^2)q^2,\\
 b_2& = (m^2 + 2mn - n^2)pq - (m^2 - 2mn - n^2)q^2,\\
c_1&= c_2 = (m^2 + n^2)p^2 + (m^2 + n^2)q^2,
\end{aligned}
\label{solabc}
\end{equation}
where $m, n, p$ and $q$ are arbitrary parameters. 

We now have two triangles with a common circumradius and a common area. For the common area to be rational, the values of $a_1, b_1, c_1$ given by \eqref{solabc} must satisfy the  condition,
\begin{equation}
(a_1+b_1+c_1)(a_1+b_1-c_1)(b_1+c_1-a_1)(c_1+a_1-b_1)=h^2, \label{areasq}
\end{equation}
where $h$ is some nonzero rational number.

On using the relations \eqref{solabc}, the condition \eqref{areasq} may be written as follows:
\begin{multline}
16mn(p^2 + q^2)^2(m + n)(m - n)(mp + mq - np + nq)(mq - np)\\
\times (mp + nq)(mp - mq + np + nq) = h^2. \label{areasq2}
\end{multline}

Now on writing, 
\begin{equation}
m=tn,\quad p=uq, \quad h=4v(u^2 + 1)n^4q^4, \label{areatr1}
\end{equation}
where $t$ is an arbitrary rational parameter, the condition \eqref{areasq2} reduces to
\begin{multline}
-(t + 1)^2(t - 1)^2t^2u^4 + (t^2 + 2t - 1)(t^2 - 2t - 1)(t - 1)(t + 1)tu^3\\
 + 6(t + 1)^2(t - 1)^2t^2u^2 - (t^2 + 2t - 1)(t^2 - 2t - 1)(t - 1)(t + 1)tu\\
 - (t + 1)^2(t - 1)^2t^2 = v^2. \label{areasq3}
\end{multline}
 
It is readily seen that when $u=1$, the left-hand side of \eqref{areasq3} is a perfect square, namely $4(t + 1)^2(t - 1)^2t^2$. We note that  the left-hand side of \eqref{areasq3} is a quartic function of $u$, and we know one value of $u$ that makes this quartic function a perfect square.  Now on applying a method described by Fermat  (as quoted by Dickson \cite[p. 639]{Di}), we  obtain    the following value of $t$ that makes the left-hand side of \eqref{areasq3} a perfect square:
\begin{multline}
u=(t^8 + 8t^7 + 20t^6 - 56t^5 - 26t^4 + 56t^3 + 20t^2 - 8t + 1)\\
\times (t^8 - 8t^7 + 20t^6 + 56t^5 - 26t^4 - 56t^3 + 20t^2 + 8t + 1)^{-1}, \label{valu}
\end{multline}

With the value of $u$ given by \eqref{valu}, and using the relations \eqref{solabc} and \eqref{areatr1}, we get the sides of two rational triangles which have a common circumradius and a common area. On appropriate scaling, the sides  $a_1, b_1, c_1$, and $a_2, b_2, c_2$ of the two triangles may be written as follows:
\begin{equation}
\begin{aligned}
a_1 &= 2t(t^4 - 2t^2 + 5)(5t^4 - 2t^2 + 1)(t^8 + 8t^7 + 20t^6 \\
& \quad \quad - 56t^5 - 26t^4 + 56t^3 + 20t^2 - 8t + 1),\\
b_1 &= (t - 1)(t + 1)(t^4 - 4t^3 + 10t^2 - 4t + 1)(t^4 + 4t^3 \\
& \quad \quad + 10t^2 + 4t + 1)(t^8 - 8t^7 + 20t^6 + 56t^5 - 26t^4\\
& \quad \quad - 56t^3 + 20t^2 + 8t + 1),\\
c_1 &= (t^2 + 1)(t^{16} + 104t^{14} - 548t^{12} + 3032t^{10} - 4922t^8\\
 & \quad \quad+ 3032t^6 - 548t^4 + 104t^2 + 1),
\end{aligned}
\end{equation}
\begin{equation}
\begin{aligned}
a_2 &= (t - 1)(t + 1)(t^4 - 4t^3 + 10t^2 - 4t + 1)(t^4 + 4t^3\\
& \quad \quad + 10t^2 + 4t + 1) (t^8 + 8t^7 + 20t^6 - 56t^5 - 26t^4\\
 & \quad \quad + 56t^3 + 20t^2 - 8t + 1),\\
b_2 &= 2t(t^4 - 2t^2 + 5)(5t^4 - 2t^2 + 1)(t^8 - 8t^7 + 20t^6 + 56t^5\\
& \quad \quad - 26t^4 - 56t^3 + 20t^2 + 8t + 1), \\
c_2 &= (t^2 + 1)(t^{16} + 104t^{14} - 548t^{12} + 3032t^{10} - 4922t^8 \\
& \quad \quad + 3032t^6 - 548t^4 + 104t^2 + 1),
\end{aligned}
\end{equation}
where $t$ is an arbitrary parameter.

The common circumradius of the two triangles is 
\begin{multline*}
(t^2+1)(t^4 - 2t^2 + 5)(5t^4 - 2t^2 + 1)(t^4 - 4t^3 + 10t^2 - 4t + 1)(t^4 + 4t^3 \\
+ 10t^2+ 4t + 1)(t^8 + 8t^7 + 20t^6 - 56t^5 - 26t^4 + 56t^3 + 20t^2 - 8t + 1)\\
\times (t^8 - 8t^7 + 20t^6 + 56t^5 - 26t^4 - 56t^3 + 20t^2 + 8t + 1)\{2(3t^4 - 6t^2 - 1)\\
\times (t^4 + 6t^2 - 3)(t^4 - 4t^3 - 6t^2 - 4t + 1)(t^4 + 4t^3 - 6t^2 + 4t + 1)\}^{-1},
\end{multline*}
while their common area is 
\begin{multline*}
t(t - 1)(t + 1)(3t^4 - 6t^2 - 1)(t^4 + 6t^2 - 3)(t^4 - 4t^3 - 6t^2 - 4t + 1)(t^4 + 4t^3 \\
- 6t^2 + 4t + 1)(t^{16} + 104t^{14} - 548t^{12} + 3032t^{10} - 4922t^8 + 3032t^6 - 548t^4 + 104t^2 + 1).
\end{multline*}

As a numerical example when $t=2$, we get two triangles with sides 
\[3283540, 7603539, 7776485, \quad {\rm and} \quad  4279155, 5834452,  7776485,\]
having common circumradius $10402718520025/2639802$ and common area $12317028393582$.

We note that additional solutions of Eq.~\eqref{areasq3} may be found by repeated application of the aforementioned method of Fermat. Alternatively, we may consider Eq.~\eqref{areasq3} as a quartic model of an elliptic curve, reduce it by a birational transformation to the cubic Weierstrass  model using the known rational point on the quartic curve \eqref{areasq3}, and find additional rational points on the cubic elliptic curve using the group law. These rational points may then be used to obtain additional  parametric solutions of our problem.

\section{Some open problems}
We have obtained pairs of rational triangles with a common  circumradius and also having either a common perimeter or a common inradius or a common area. It would be of interest to determine whether there are three or more rational triangles with  a common  circumradius and  also having  a common perimeter or a common inradius or a common area.

\noindent Ajai Choudhry, 13/4 A Clay Square, Lucknow - 226001, India.

\noindent E-mail address: ajaic203@yahoo.com

\end{document}